\documentclass[10pt,reqno]{amsart}
\usepackage[utf8]{inputenc}
\usepackage[english]{babel}

\usepackage{amsmath, amsfonts, amsthm, amssymb, stmaryrd, color, cite}
\usepackage{parskip}

\usepackage{latexsym,hyperref}

\usepackage{cite}
\usepackage{amsmath, amsfonts, amsthm, amssymb, stmaryrd, color, cite}

\usepackage{amsthm}
\usepackage{cite}
\usepackage{color}
\usepackage[dvipsnames]{xcolor}

\author{}
\newtheorem{theorem}{Theorem}[section]

\title[\emph{Some Exit Time Estimates for SBM and FVP}]
{Some Exit Time Estimates for Super-Brownian Motion and Fleming-Viot Process}

\author[P. Fatheddin]{Parisa Fatheddin}
\address{Department of Mathematics, University of Pittsburgh, Pittsburgh, PA 15260, USA.}
\email{PAF49@pitt.edu}


\begin{document}
\maketitle
\begin{abstract}
Estimates for exit time from an interval of length $2r$ before a prescribed time $T$ are derived for solutions of a class of stochastic partial differential equations used to characterize two population models: super-Brownian motion and Fleming-Viot Process. These types of estimates are then derived for the two population models. The corresponding large deviation results are also applied for the acquired bounds.
\end{abstract}

\section{Introduction}
\setcounter{equation}{0}

\indent Since the early works of M. Freidlin and A. Wentzell \cite{Freidlin1, Freidlin3}, many authors have investigated the exit problem from a bounded domain. These estimates have offered notable insight in the fields of applications, where exit times can be translated to determining the probability of a particular quantity exceeding a specific threshold. In finance, authors in \cite{Dshalalow} use their results on exit time to determine the time that the stock price exceeds a particular price and if an option is not exercised by a specific deadline. In communication theory, for example in \cite{Cottrell}, a radio channel is considered in which messages are transmitted between users and the exit time is given to determine when the number of blocked users reach a critical value that would break the system. Furthermore, for applications in queueing theory, we refer the reader to \cite{Dupuis} and Chapter 11 of \cite{Shwartz}. The importance of exit problem in population models is intuitively clear. Here we consider the exit measure for well studied population models in the literature: super-Brownian motion and Fleming-Viot process. First we consider a class of stochastic partial differential equations used in \cite{Xiong} to characterize the two models and give its solution's exit time. Since the results in \cite{Xiong} were achieved only in dimension one, our exit problem is limited to this case only.\\
In \cite{Dynkin} exit measure was studied for superprocesses, of which SBM and FVP are two cases. As for the exit measure of SBM, based on ideas provided in \cite{Dynkin}, authors in \cite{Hesse, Sheu, Verzani} considered an increasing sequence $D_{k}$ of subdomains in the bounded domain of study, $D$ and studied the probability of the process hitting these subdomains. Authors in \cite{Le Gall} proved that only in dimensions one and two, boundary points of $D$ get hit with positive probability and this probability is zero for higher dimensions.\\
Moreover, for FVP, authors in \cite{Bieniek, Grigorescu, Burdzy} considered the situation in which particles are destroyed upon hitting the boundary of the domain and since the size of the population in FVP is assumed to be constant, then at the occurrence of this event, another particle in the domain reproduces one offspring. Investigators in \cite{Grigorescu} further studied the control of particles hitting the boundary and established the tightness of the average number of visits to the boundary. To the best of our knowledge, estimates for the time of exit measure of FVP have not been previously shown in the literature.\\
Also it is an interesting problem to consider exit time if there exists an attraction point in the domain. In the context of populations, this attraction point can be interpreted as a food source. Mathematically, an attraction point is formulated as a point to which trajectories of the solution converge to after sufficiently long period of time (i.e. $t\rightarrow \infty$). Following ideas from Chapter 4 of \cite{Freidlin1} we derive an upperbound for the probability of the first time the trajectory enters a $\delta$-neighborhood of the attraction point after a prescribed time, and also an upperbound for the mean exit time from the domain in the case of the existence of an attraction point. We determine these results for the class of SPDEs and the two population models.\\
The study of exit time emerged from the theory of large deviations and one may obtain estimates on exit time using the corresponding large deviation principle as is described in Section 5.7 of \cite{Dembo} and performed in \cite{Chenal, Cerrai, Jung}. Here we use this connection to derive our estimates for the exit problem of the solutions of class of SPDEs and population models. As mentioned above, a few other authors have also examined the exit measure of super-Brownian motion (SBM) and Fleming-Viot Process (FVP). For both models, our method and results are new to the literature, since they rely on a direct approach based on the characterization of the models by the unique solution to stochastic partial differential equations introduced in \cite{Xiong}. In addition, to the best of our knowledge, the connection between exit measure of the two population models with their corresponding large deviations results have not previously been given in the literature.\\
This article is organized as follows. We provide a brief background to the population models studied here along with the statement of main results in Section 2. Then in Section 3, estimates for exit times are derived for the class of SPDEs using their connections with large deviations. These estimates are afterwards achieved for SBM and FVP in Section 4.

\section{Notation and Main Results}

We begin by giving a brief description on the two population models considered here. Super-Brownian motion (SBM) is the continuous version of branching Brownian motion where population evolves as a cloud and branches out like a tree. With regards to applications in biology and sociology, super-Brownian motion has been used to investigate Brownian bug model and voter model as in \cite{Bramson, Birch}. Brownian bug model studies small organisms such as bugs that reproduce by binary fission at a constant rate and organisms move according to a random walk. When it is assumed that the birth and death of bugs are independent of the spatial distribution of the population and bugs are assumed to move according to a Brownian motion, then the model becomes a super-Brownian motion, which in this context is also referred to as the Felsenstein's problem. Voter model assumes that each individual in the population has an opinion a or b and it examines the spread of one opinion over the other. It has been shown in \cite{Cox} that voter model converges weakly to super-Brownian motion.\\
   As for the other population model studied here, Fleming-Viot process (FVP) is the continuous version of step-wise mutation model, in which individuals are grouped based on their gene type. This process initially was developed in the study of diffusion models in population genetics in the paper \cite{Fleming} by Fleming and Viot in 1979. Many authors have since considered this process to study different biological developments for example, alleles diffusions in \cite{Ethier} and polarity in cells in \cite{Gupta}, which is the tendency of the majority of membrane molecules in a cell to concentrate in one place, a condition that occurs mostly in yeast cells. Another important application in biology is the modelling of parasite metapopulations by Fleming-Viot process as in \cite{Taylor}, where parasites are assumed to spread and infect their host like the spread of a disease. For more applications in population genetics we recommend \cite{Ethier}.
      The distribution of gene types is studied in FVP making it a probability measure-valued process; whereas, SBM is a measure-valued process. For more information and background on these population models we refer the reader to \cite{Dawson, Etheridge, moderate, Perkins}.\\
We now give the stochastic partial differential equation (SPDE), the existence and uniqueness of solutions of which were studied in \cite{Xiong} in dimension one and was used to characterize SBM and FVP as unique solutions to certain SPDEs. Suppose $(\Omega, \mathcal{F}, P)$ is a probability space with $\{\mathcal{F}_{t}\}$ being a family of non-decreasing sub $\sigma$-fields with the standard conditions of $\{\mathcal{F}_{t}\}$ being right continuous and $\mathcal{F}_{0}$ $P$-complete. Since super-Brownian motion is a measure-valued process, we may denote it as $\{\mu_{t}^{\epsilon}\}$, where $\epsilon$ is the branching rate. Using $u_{t}^{\epsilon}(y)= \int_{0}^{y}\mu_{t}^{\epsilon}(dx)$ for $y\in \mathbb{R}$, this population model was given in \cite{Xiong} by,
\begin{equation}\label{SBM}
u_{t}^{\epsilon}(y)= F(y)+ \sqrt{\epsilon}\int_{0}^{t}\int_{0}^{u_{s}^{\epsilon}(y)}W(dads) + \int_{0}^{t}\frac{1}{2}\Delta u_{s}^{\epsilon}(y)ds,
\end{equation}
where $F(y)= \int_{0}^{y}\mu_{0}(dx)$ and $W$ is an $\mathcal{F}_{t}$-adapted space-time white noise random measure on $\mathbb{R}^{+} \times \mathbb{R}$ with intensity measure $dsda$. Similarly, Fleming-Viot Process, $\{\mu_{t}^{\epsilon}\}$, where $\epsilon$ is the mutation rate, is a probability measure-valued process, and was characterized by
\begin{equation}\label{FVP}
u_{t}^{\epsilon}(y)= F(y) + \sqrt{\epsilon} \int_{0}^{t} \int_{0}^{1} \left(1_{a\leq u_{s}^{\epsilon}(y)}-u_{s}^{\epsilon}(y)\right) W(dads) + \int_{0}^{t} \frac{1}{2} \Delta u_{s}^{\epsilon}(y)ds,
\end{equation}
using, $u_{t}^{\epsilon}(y)= \mu_{t}^{\epsilon}((-\infty, y])$ with $F$ a function on $\mathbb{R}$ and $W$ a space-time white noise random measure on $\mathbb{R}^{+} \times [0,1]$. In \cite{Xiong}, observing the similarities in equations \eqref{SBM} and \eqref{FVP}, the following equation,
\begin{equation}\label{SPDE}
u_{t}^{\epsilon}(y)= F(y)+ \sqrt{\epsilon}\int_{0}^{t}\int_{U}G(a,y,u_{s}^{\epsilon}(y))W(dads) + \int_{0}^{t}\frac{1}{2}\Delta u_{s}^{\epsilon}(y)ds,
\end{equation}
was introduced, where $W$ is an $\mathcal{F}_{t}$-adapted space-time white noise random measure on $\mathbb{R}^{+}\times U$, $F$ is a function on $\mathbb{R}$ and $(U,\mathcal{U},\lambda)$ is a measure space with $\lambda$ denoting the Lebesgue measure with $U=\mathbb{R}$ for SBM and $U=[0,1]$ in the case of FVP. Furthermore, $G:U\times \mathbb{R}^{2}\rightarrow \mathbb{R}$ is H\"older continuous and satisfies the linear growth condition:
\begin{eqnarray}
\int_{U}\left|G(a,y,u_{1})-G(a,y,u_{2})\right|^{2}\lambda(da) &\leq& K|u_{1}-u_{2}|, \label{condition1}\\
\int_{U}\left|G(a,y,u)\right|^{2}\lambda(da)&\leq& K\left(1+|u|^{2}\right). \label{condition2}
\end{eqnarray}For $0<\beta \in \mathbb{R}$, let $\mathcal{M}_{\beta}(\mathbb{R})$ be the set of $\sigma$-finite measures $\mu$ on $\mathbb{R}$ such that
\begin{equation*}
\int e^{-\beta |x|} d\mu(x)<\infty.
\end{equation*}
Let $\alpha \in (0,1)$ and denote $\mathbb{B}_{\alpha,\beta}$ to be the space of all functions $f:\mathbb{R}\rightarrow \mathbb{R}$ such that for all $m\in \mathbb{N}$,
\begin{eqnarray}
\left|f(y_{1})-f(y_{2})\right|&\leq& Ke^{\beta m}|y_{1}-y_{2}|^{\alpha} \hspace{.3cm} \text{ $\forall |y_{1}|,|y_{2}|\leq m$,}\\ \label{con1}
|f(y)|&\leq& Ke^{\beta |y|} \hspace{.3cm} \text{ $\forall y\in \mathbb{R}$}. \label{con2}
\end{eqnarray}
 In \cite{large} it was assumed that $F(y)$ is in space $\mathbb{B}_{\alpha,\beta_{0}}$ for $\alpha \in (0,\frac{1}{2})$, $\beta_{0}\in (0,\beta)$ and the class of SPDEs given by \eqref{SPDE} was proved to be in space $\mathcal{C}\left([0,1];\mathbb{B}_{\beta}\right)$, where $\mathbb{B}_{\beta}$ is the collection of continuous functions on $\mathbb{R}$ satisfying only condition \eqref{con2}. Furthermore, it can be shown that $\mathbb{B}_{\beta}$ is a Banach space with norm,
\begin{equation}\label{norm1}
\|f\|_{\beta} = \sup_{y\in \mathbb{R}} e^{-\beta |y|}|f(y)|.
\end{equation}
Since the existence and uniqueness of solutions achieved in \cite{Xiong} are limited to one dimension, then we consider the domain of study to be the interval $(-r,r)$ and the first exit time to be denoted as $\tau^{\epsilon}:= \inf\{t: u_{t}^{\epsilon}(y)\notin D\}$. We derive estimates for exit times directly then use our results on large deviations with rate function, $I_{1}: \tilde{\mathcal{C}}_{\beta} \rightarrow [0,\infty]$ given later in the section. For better presentation, we denote,
\begin{equation*}
J(r,\epsilon, T):= \sup_{0<t\leq T} \frac{8k \epsilon^{\frac{1}{2}} C_{2}^{\frac{1}{k}}C_{3}\sqrt{T}}{(r\sqrt{t}- C_{1}C_{4})(k-1)},
\end{equation*}
for positive constants, $C_{1},C_{2},C_{3},$ and $C_{4}$. Using this notation, we have the following estimates.
\begin{theorem}\label{exitSPDE}
For the family, $\{u_{t}^{\epsilon}(y)\}_{\epsilon>0}$ given by \eqref{SPDE}, and a given $\delta >0$, there exists $\epsilon_{0} >0$ such that for all $0<\epsilon <\epsilon_{0}$,\\
a. the likelihood of the exit time of a trajectory from the domain $(-r,r)$ to be before a given time $T$ is approximated by,
\begin{equation}\label{probabilitySPDE}
  \exp\left(-\frac{1}{\epsilon}\left(\inf_{u_{t}(y)\in [-r,r]^{c}}I_{1}(u_{t}(y)) - \delta\right)\right)\leq P(\tau^{\epsilon}\leq T) \leq  J(r,\epsilon, T),
\end{equation}
b. for mean exit time we have,
\begin{equation}\label{expectationSPDE}
 \mathbb{E}(\tau^{\epsilon}) \leq \frac{1}{1-\exp \left( -\frac{1}{\epsilon} \left(\inf_{u_{t}(y) \in [-r,r]} I_{1}\left(u_{t}(y)\right)- \delta\right)\right)}.
\end{equation}
\end{theorem}

\begin{theorem}\label{attractionpoint}
Suppose domain $(-r,r)$ possesses exactly one attraction point being zero, then for a given $\delta>0$ and all $\epsilon >0$ sufficiently small, \\
a. the likelihood of the first exit time being before a prescribed time, $T$ is estimated by,
\begin{equation*}
  1-J(\delta_{0},\epsilon, T) \leq P(\tau^{\epsilon} \leq T)\leq J(r,\epsilon, T) + \exp\left( -\frac{1}{\epsilon}\left(\inf_{u_{t}(y)\in (-\delta_{0},\delta_{0})^{c}} I_{1}(u_{t}(y)) -\delta\right)\right),
\end{equation*}
b. the mean exit time is bounded by,
\begin{equation*}
\mathbb{E}(\tau^{\epsilon}) \leq \frac{1}{1- \exp\left(-\frac{1}{\epsilon}\left(\inf_{u_{t}(y) \in (-\delta_{0},\delta_{0})^{c}}I_{1}(u_{t}(y)) - \delta\right)\right)}.
\end{equation*}
\end{theorem}

For the next theorems regarding the two populations models under study, we define,
\begin{equation}\label{innerproduct}
<f,\mu_{t}^{\epsilon}>_{\beta} := \int \sup_{y} e^{-\beta |y|} |f(y)|\mu_{t}^{\epsilon}(dy),
\end{equation}
and let
\begin{equation}\label{N2}
 N_{2} := \frac{1}{K_{9}\sqrt{\epsilon T} (1+e^{3\beta_{1}|r|})}\left(1- \sup_{0<t\leq T} \frac{1}{t} \left(K_{5} + K_{6} e^{2\beta_{0}|r|}\right)\right),
 \end{equation}
for positive constants $K_{5}, K_{6}$ and $K_{9}$.

\begin{theorem}\label{Meansize}
For both cases of super-Brownian motion and Fleming-Viot Process, the exit time is approximated by,
\begin{equation}\label{result}
P\left(\tau^{\epsilon} \leq T\right) \leq \sup_{0<t\leq T} \frac{8k\epsilon^{\frac{1}{2}} C_{2}^{\frac{1}{k}}\sqrt{T}}{(\sqrt{N_{2}t} -C_{1}C_{5})(k-1)},
\end{equation}
where $C_{5}:= \sup_{y\in \mathbb{R}} e^{-(\beta_{1}-\beta_{0})|y|}$. Moreover, the mean of the size of the population at time $t$ for both population models is estimated by,
\begin{eqnarray}\label{result2}
\mathbb{E}\|\mu_{t}^{\epsilon}\|_{\beta}^{2} &\leq& M(t^{2}+t^{3}) \left(K\beta^{6} + K\beta^{4} - K\beta^{2}\right) + \ln t \nonumber\\
&&+ K(\beta -\beta_{0})^{2}(t+1) + KM\epsilon \left(\sqrt{t} + t^{3/2}\right).
\end{eqnarray}
\end{theorem}

 For completeness, we provide the results on large deviations achieved in \cite{large} as follows. Because of the nonlipschitz continuity nature of the coefficient $G(a,y,u_{s}^{\epsilon}(x))$, the existence and uniqueness of solutions could not be obtained for the controlled PDE, also referred to as the skeleton equation given below,
\begin{equation}\label{controlled}
u_{t}(y)= F(y) + \int_{0}^{t}\int_{U}G(a,y,u_{s}(y))h_{s}(a)\lambda(da)ds + \int_{0}^{t} \frac{1}{2} \Delta u_{s}(y)ds,
\end{equation}
where $h_{s}(\cdot)\in L^{2}\left([0,1]\times U, ds\lambda(da)\right)$. As a consequence, an equivalence class was introduced to group solutions in a suitable way. That is, let $u\sim_{1} v$ if both $u$ and $v$ are solutions to \eqref{SPDE} with the same function $h_{s}(\cdot)$. Then $\mathcal{C}\left([0,1];\mathbb{B}_{\beta}\right)/\sim_{1}$ is a pseudo-metric space and to convert this space to a Polish space, another equivalence class was applied defined by $x\sim_{2} y$ if $d(x,y)=0$. Namely, large deviation principle was achieved for the class of SPDEs in
$\tilde{\mathcal{C}}_{\beta}:=\bar{\mathcal{C}}\left([0,1];\mathbb{B}_{\beta}\right)/ \sim_{1} / \sim_{2}$. For SPDE \eqref{SPDE} this result was given as Theorem 2 in \cite{large} in space $\tilde{\mathcal{C}}_{\beta}$ with rate function $I_{1}(\hat{u})= \inf_{u\in \hat{u}} I(u)$ where,
  \begin{equation}\label{2.7}
I(u) =\left\{\begin{array}{ll}  \frac{1}{2} \inf\left\{ \displaystyle
\int_{0}^{1} \int_{U} |h_{s}(a)|^{2} \lambda(da)ds : u = \gamma (h)\right\}
 &\exists h\; \mbox{s.t. } u=\gamma(h),\\
\infty &\mbox{otherwise.}
\end{array}\right.
\end{equation}
Here $\gamma$ is a map from $L^{2}([0,1] \times U, ds\lambda(da))$ to $\tilde{\mathcal{C}}_{\beta}$ with domain consisting of $h$ such that \eqref{controlled} has a solution, and the equivalence class of the solution is denoted as $u=\gamma(h)$. As for SBM and FVP, Cameron-Martin space, $H_{\nu}$ was introduced to give large deviations in spaces $\tilde{\mathbb{M}}_{\beta}:=\bar{\mathcal{C}}\left([0,1];\mathcal{M}_{\beta}(\mathbb{R})\right)/ \sim_{1} / \sim_{2}$ and $\tilde{\mathbb{P}}_{\beta}:=\bar{\mathcal{C}}\left([0,1];\mathcal{P}_{\beta}(\mathbb{R})\right)/ \sim_{1} / \sim_{2}$, respectively. We provide the conditions for Cameron-Martin space as follows. For $\nu$ on the set of finite measures on $\mathbb{R}$, $\mathcal{M}_{F}(\mathbb{R})$, the Cameron-Martin space, $H_{\nu}$, is the set of measures $\mu \in C([0,1];\mathcal{M}_{F}(\mathbb{R}))$ such that,
 \begin{enumerate}
\item  $\mu_{0}=\nu$,
\item the $\mathcal{D}^{*}$-valued map $t\mapsto \mu_{t}$ defined on [0,1] is absolutely continuous with
 respect to time. Let $\dot{\mu}$ and $\Delta^{*}\mu$ be its generalized derivative and Laplacian
  respectively,
\item for every $t\in [0,1]$, $\dot{\mu}_t - \frac{1}{2}\Delta^{*} \mu_{t} \in \mathcal{D}^{*}$ is
 absolutely continuous with respect to $\mu_{t}$
 with $\frac{d\left(\dot{\mu}_{t}-\frac{1}{2}\Delta^{*}\mu_{t}\right)}{d\mu_t}$
  in $L^{2}([0,1] \times \mathbb{R},
ds\mu(dy))$ being the (generalized) Radon Nikodym derivative.
\end{enumerate}
In Theorem 3 of \cite{large}, large deviations was achieved for SBM with rate function, $I_{2}(\hat{\mu})= \inf_{\mu \in \hat{\mu}} I(\mu)$ where,
 \begin{equation}\label{rate4sbm}
  I(\mu)=  \left\{\begin{array} {ll}  \frac{1}{2} \displaystyle \int_{0}^{1}
  \int_{\mathbb{R}}\left|\frac{\left(\dot{\mu}_{t} - \frac{1}{2} \Delta^* \mu_{t}\right)(dy)}{\mu_{t}(dy)}\right|^2 \mu_{t}(dy) dt
   & \mbox{\emph{if }} \mu \in H_{\mu_{0}}, \\
   \infty
  & \mbox{\emph{otherwise}.}     \end{array}   \right.
 \end{equation}
Similarly in Theorem 4 of \cite{large}, large deviation principle was given for FVP with rate function, $I_{3}(\hat{\mu})= \inf_{\mu \in \hat{\mu}} I(\mu)$ where,
 \begin{equation}\label{rate4fvp}
  I(\mu)=  \left\{\begin{array} {ll}  \frac{1}{2} \displaystyle \int_{0}^{1}
  \int_{\mathbb{R}}\left|\frac{\left(\dot{\mu}_{t} - \frac{1}{2} \Delta^* \mu_{t}\right)(dy)}{\mu_{t}(dy)}\right|^2 \mu_{t}(dy) dt
   & \mbox{\emph{if }} \mu \in \tilde{H}_{\mu_{0}}, \\
   \infty
  & \mbox{\emph{otherwise}.}     \end{array}   \right.
 \end{equation}
where $\tilde{H}_{\nu}$ is the Cameron Martin space with conditions involving probability measures $\mathcal{P}(\mathbb{R})$ instead of $\mathcal{M}_{F}(\mathbb{R})$ and the extra condition,
\begin{equation*}
\left<\mu_{t}, \frac{d\left(\dot{\mu}_{t}-\frac{1}{2}\Delta^{*}\mu_{t}\right)}{d\mu_{t}}\right>=0.
\end{equation*}
For each of these results $F(y)$ is assumed to be in space $\mathbb{B}_{\alpha,\beta_{0}}$ for $\alpha \in (0,\frac{1}{2})$ and $\beta_{0} <\beta$.

\section{Exit Times for the Class of SPDEs}

In this section, we derive estimates for exit times for the class of SPDEs that was introduced in the previous section. Namely, our aim here is to prove Theorems \ref{exitSPDE} and \ref{attractionpoint}. Since we have the uniqueness of strong solutions to SPDE \eqref{SPDE} established in \cite{Xiong}, we may use its mild form,
\begin{equation}\label{mild}
u_{t}^{\epsilon}(y)= \int_{\mathbb{R}}p_{t}(y-x)F(x)dx + \sqrt{\epsilon}\int_{0}^{t}\int_{U}\int_{\mathbb{R}}p_{t-s}(y-x)G(a,x,u_{s}^{\epsilon}(x))dxW(dads),
\end{equation}
with the Brownian semigroup, where $p_{t}(x)=\frac{1}{\sqrt{2\pi}t}e^{-\frac{x^{2}}{2t}}$ is the heat kernel. Let the first exit time from the domain, $(-r,r)$ be denoted as $\tau^{\epsilon}:= \inf\{t: u_{t}^{\epsilon}(y)\notin (-r,r)\}$.

\textbf{Proof of Theorem 2.1}: We begin by proving \eqref{probabilitySPDE}, where we apply a direct method for the upperbound based on ideas in the proof of Theorem 2 in \cite{Nualart} and for its lowerbound we use the large deviation result stated in Section 2. Observe that for a prescribed time $T$,
\begin{equation}\label{tau}
P(\tau^{\epsilon} >T) = P\left(\sup_{0 < t\leq T} \|u_{t}^{\epsilon}(y)\|_{\beta} < r\right),
\end{equation}
where we have excluded $t=0$ since initially the solution is assumed to be in the domain. Using the $\beta$ norm given by \eqref{norm1}, we proceed to find,
\begin{eqnarray}\label{estimate}
&&P\left(\sup_{0 < t\leq T} \|u_{t}^{\epsilon}(y)\|_{\beta} \geq r\right)\\
&\leq& P\left( \sup_{0 < t\leq T} \sup_{y\in \mathbb{R}} e^{-\beta |y|} \left|\int_{\mathbb{R}} p_{t}(y-x) F(x)dx\right|\right. \nonumber\\
&&\left.+  \sup_{0 < t\leq T} \sup_{y\in \mathbb{R}} e^{-\beta |y|} \left|\sqrt{\epsilon} \int_{0}^{t} \int_{U} \int_{\mathbb{R}} p_{t-s}(y-x) G(a,x,u_{s}^{\epsilon}(x))dx W(dads)\right|\geq r \right). \nonumber
\end{eqnarray}
Recall that $F$ is assumed to be in $\mathbb{B}_{\alpha,\beta_{0}}$ space with $\alpha \in (0,\frac{1}{2})$ and $\beta_{0}\in (0,\beta)$ satisfying inequality \eqref{con2} with $\beta$ replaced by $\beta_{0}$ and constant denoted here as $K_{1}$. Thus,
\begin{eqnarray}\label{Finequality}
\sup_{y\in \mathbb{R}} e^{-\beta |y|} \left|\int_{\mathbb{R}} p_{t}(y-x) F(x)dx\right| &\leq& \sup_{y\in \mathbb{R}} e^{-\beta |y|}\int_{\mathbb{R}} p_{t}(y-x) K_{1} e^{\beta_{0}|x|}dx \nonumber\\
&\leq& \frac{K_{1}K_{2}}{\sqrt{2\pi t}} \sup_{y\in \mathbb{R}} e^{-(\beta-\beta_{0}) |y|},
\end{eqnarray}
where
\begin{equation*}
K_{2}:= \int_{\mathbb{R}} \exp\left(-\frac{(y-x)^{2}}{2t} - (|y|-|x|)\beta_{0}\right)dx.
\end{equation*}
Let $I(t)= \sqrt{\epsilon}\int_{0}^{t} \int_{U} \int_{\mathbb{R}} p_{t-s}(y-x) G(a,x,u_{s}^{\epsilon}(x))dxW(dads)$, then letting $C_{1}= \frac{K_{1}K_{2}}{\sqrt{2\pi}}$, the probability given by \eqref{estimate} becomes bounded by,
\begin{equation}\label{expression}
P\left(\sup_{0 < t\leq T} \sup_{y\in \mathbb{R}} e^{-\beta |y|} |I(t)| \geq r - \sup_{0 < t\leq T} \frac{C_{1}}{\sqrt{t}} \sup_{y\in \mathbb{R}} e^{-(\beta-\beta_{0}) |y|} \right).
\end{equation}
For consecutive approximations we need the following estimate which was established in \cite{large} as Lemma 1. For any $n\geq 1$, and $\beta_{1} \in \left(\beta_{0},\beta\right)$,
\begin{equation}\label{Mlemma}
M:= \mathbb{E}\left(\sup_{0\leq s\leq 1} \int_{\mathbb{R}} |u_{s}^{\epsilon}(x)|^{2}e^{-2\beta_{1}|x|}dx\right)^{n} <\infty.
\end{equation}
Inspired by the proof of Theorem 2 in \cite{Nualart} and similarly Proposition 7 in \cite{Rovira}, we apply the well-known Kolmogorov criterion (see for example Corollary 1.2 in \cite{Walsh}). Denoting,
\begin{equation*}
P(t,t_{1}):= p_{t-s}(y-x)-p_{t_{1}-s}(y-x),
\end{equation*}
we have for $k\geq 2$, and $t_{1}\leq t$,
\begin{eqnarray*}
\mathbb{E}\left|I(t)-I(t_{1})\right|^{k} &\leq& 2^{k-1} \mathbb{E} \left|\sqrt{\epsilon} \int_{0}^{t} \int_{U} \int_{\mathbb{R}} P(t,t_{1}) G(a,x,u_{s}^{\epsilon}(x))dxW(dads)\right|^{k}\\
&&+ 2^{k-1} \mathbb{E} \left|\sqrt{\epsilon} \int_{t_{1}}^{t} \int_{U} \int_{\mathbb{R}} p_{t_{1}-s}(y-x) G(a,x,u_{s}^{\epsilon}(x)) dx W(dsda) \right|^{k} \\
&=& J_{1} + J_{2}.
\end{eqnarray*}
Using the Burkholder-Davis-Gundy and H\"older inequalities along with inequality \eqref{condition2} we arrive at,
\begin{eqnarray*}
J_{1} &\leq& 2^{k-1} \epsilon^{\frac{k}{2}} \mathbb{E} \left|\int_{0}^{T} \int_{U} \left(\int_{\mathbb{R}} 1_{[0,t]} (s) P(t,t_{1}) G(a,x,u_{s}^{\epsilon}(x))dx\right)^{2} dads\right|^{\frac{k}{2}}\\
&\leq& 2^{k-1} \epsilon^{\frac{k}{2}} \mathbb{E} \left|\int_{0}^{t} \int_{U}\int_{\mathbb{R}} P(t,t_{1})^{2}e^{2\beta_{1} |x|} dx \int_{\mathbb{R}} G(a,x,u_{s}^{\epsilon}(x))^{2}e^{-2\beta_{1} |x|}dx dads \right|^{\frac{k}{2}}\\
&\leq& 2^{k-1} \epsilon^{\frac{k}{2}} \mathbb{E} \left| \int_{0}^{t} \int_{\mathbb{R}} P(t,t_{1})^{2} e^{2\beta_{1} |x|}dx \int_{\mathbb{R}} \left(1+ |u_{s}^{\epsilon}(x)|^{2}\right)e^{-2\beta_{1} |x|}dx ds\right|^{\frac{k}{2}}.
\end{eqnarray*}
In the proof of Lemma 4 in \cite{large}, it was found that for $\alpha \in (0,\frac{1}{2})$,
\begin{equation}\label{Pbound}
\tilde{J}_{1}(s):= \int_{\mathbb{R}} P(t_{1},t_{2})^{2} e^{2\beta_{1} |x|} dx \leq K e^{2 \beta_{1}|y|} |t_{1}-t_{2}|^{\alpha},
\end{equation}
where the authors wrote $P(t_{1},t_{2})$ as,
\begin{equation*}
P(t_{1},t_{2})^{2}= \left|p_{t_{1}-s} (y-x) - p_{t_{2}-s}(y-x)\right|^{\alpha} \left|p_{t_{1}-s}(y-x)-p_{t_{2}-s}(y-x)\right|^{2-\alpha}.
\end{equation*}
Interval $(0, \frac{1}{2})$ was used there to obtain the suitable bound. We note that each estimate performed in achieving \eqref{Pbound} also holds for $\alpha =1$. Therefore,
\begin{equation*}
J_{1}\leq  2^{k-1} \epsilon^{\frac{k}{2}}MK_{3} e^{\beta_{1}|y|k}|t-t_{1}|^{\frac{k}{2}},
\end{equation*}
where constant $M$ is given by \eqref{Mlemma}. Similarly,
\begin{eqnarray*}
J_{2} &\leq& 2^{k-1} \epsilon^{\frac{k}{2}}\mathbb{E} \left|\int_{t_{1}}^{t}  \int_{\mathbb{R}} p_{t_{1}-s}(y-x)^{2} e^{2\beta_{1} |x|}dx \int_{\mathbb{R}} \left(1+ |u_{s}^{\epsilon}(x)|^{2}\right)e^{-2\beta_{1} |x|}dx\right|^{\frac{k}{2}} \\
&\leq& 2^{k-1}\epsilon^{\frac{k}{2}}MK_{4} e^{\beta_{1}|y|k} |t-t_{1}|^{\frac{k}{2}}.
\end{eqnarray*}
Leading to,
\begin{equation*}
\mathbb{E}|I(t)-I(t_{1})|^{k}\leq 2^{k-1}  \epsilon^{\frac{k}{2}} M e^{\beta_{1}|y|k} (K_{3}+K_{4}) |t-t_{1}|^{\frac{k}{2}}.
\end{equation*}
Next we let $\psi(x)=|x|^{k}$ for $k\geq 2$, and $p(x)= e^{\beta_{1}|y|}x^{2}$ for $x\geq 0$ and denote,
\begin{equation*}
B:= \int_{0}^{T} \int_{0}^{T} \frac{ |I(t)-I(t_{1})|^{k}}{p(|t-t_{1}|)^{k}} dt_{1}dt,
\end{equation*}
then as was shown above,
\begin{eqnarray*}
\mathbb{E}(B) &\leq&  (K_{3}+K_{4})2^{k-1} \epsilon^{\frac{k}{2}} M \int_{0}^{T}\int_{0}^{T} |t-t_{1}|^{\frac{k}{2}-2k} dt_{1}dt\\
&=&    \left\{\begin{array} {ll}
\frac{(K_{3}+K_{4})2^{k+2} \epsilon^{\frac{k}{2}} M }{9k^{2}-18k+8} T^{\frac{-3k+4}{2}}
   & \mbox{\emph{if }} \frac{-3k+4}{2} \mbox{ is even} \\
   0
  & \mbox{\emph{if }} \frac{-3k+4}{2} \mbox{ is odd}    \end{array}   \right.
\end{eqnarray*}
Since the preceding estimates hold for any $k\geq 2$, then we may assume that $k\geq 2$ is such that $\frac{-3k+4}{2}$ is even to obtain,
\begin{equation*}
\mathbb{E}(B)\leq C_{2} \epsilon^{\frac{k}{2}} T^{\frac{4-3k}{2}}.
\end{equation*}
 Furthermore, functions $\psi(x)$ and $p(x)$ are symmetric about $0$ with $\psi(x)$ being convex and $\lim_{x\rightarrow \infty} \psi(x)=\infty$. Also $p(0)=0$ with $p(x)$ increasing for $x>0$. Therefore, according to Garsia, Rodemich, Rumsey Lemma given as Theorem 1.1 in \cite{Walsh}, we have,
\begin{equation*}
|I(t)-I(t_{1})|\leq 8 \int_{0}^{|t-t_{1}|} \frac{B^{\frac{1}{k}}}{x^{\frac{2}{k}}}dp(x),
\end{equation*}
with $n=1$ since we have dimension one. We note that this lemma is in deterministic setting;however, as commented in \cite{Walsh} it is also applicable when $\mathbb{E}(B)<\infty$. Since the above holds for any $0\leq t_{1} \leq t$, we let $t_{1}=0$ to obtain,
\begin{equation*}
\sup_{0 < t\leq T} |I(t)|\leq 8B^{\frac{1}{k}} \int_{0}^{T} x^{-\frac{2}{k}}dp(x)\leq
\frac{8k}{k-1} e^{\beta_{1}|y|} T^{\frac{2k-2}{k}} B^{\frac{1}{k}}.
\end{equation*}
We denote $C_{3}:= \sup_{y\in \mathbb{R}} e^{-(\beta - \beta_{1})|y|}$ and $C_{4}:=\sup_{y\in \mathbb{R}} e^{-(\beta - \beta_{0})|y|}$ then using \eqref{expression} and noting that $B$ is a positive random variable, we apply the Markov inequality as follows,
\begin{eqnarray}\label{probability}
P\left(\sup_{0 < t\leq T}\|u_{t}^{\epsilon}(y)\|_{\beta} \geq  r\right)&\leq&
P\left( C_{3} \frac{8k}{k-1} T^{\frac{2k-2}{k}} B^{\frac{1}{k}} \geq r - \sup_{0 < t\leq T} \frac{C_{1}}{\sqrt{t}} C_{4}\right) \nonumber \\
&\leq& P\left( B^{\frac{1}{k}} \geq \frac{ r(k-1)}{8kC_{3} } T^{\frac{2-2k}{k}} - \sup_{0 < t\leq T} \frac{C_{1} C_{4}(k-1)}{8kC_{3}\sqrt{t}} T^{\frac{2-2k}{k}}\right) \nonumber \\
&\leq& \mathbb{E}(B^{\frac{1}{k}}) \sup_{0<t\leq T} \frac{8kC_{3}\sqrt{t}} {(r\sqrt{t}- C_{1}C_{4})(k-1)T^{\frac{2-2k}{k}}} \nonumber \\
&\leq&  \sup_{0<t\leq T} \frac{8k \epsilon^{\frac{1}{2}} C_{2}^{\frac{1}{k}}C_{3}\sqrt{T}}{(r\sqrt{t}- C_{1}C_{4})(k-1)},
\end{eqnarray}
where we have used the concave property of $\psi^{-1}(x)$. This estimate confirms the fact that as the given radius of exit domain becomes sufficiently large or as the noise goes to zero by setting $\epsilon \rightarrow 0$, one does not expect the solution to exit the prescribed domain. Moreover, \eqref{probability} verifies mathematically that if the prescribed time $T$ is zero, then the probability of exiting the domain before $T$ is also zero. Now denoting the right hand side of inequality \eqref{probability} as $J(r,\epsilon, T)$ and recalling \eqref{tau}, we obtain,
\begin{equation*}
P(\tau^{\epsilon} \leq T) = P\left(\sup_{0<t\leq T} \|u_{t}^{\epsilon}(y)\|_{\beta} \geq r\right) \leq J(r,\epsilon, T).
\end{equation*}
This offers an upperbound on the probability that the solution will not exit before time $t=T$.\\
Relating to the large deviation result provided in Section 2 we achieve the upperbound employing the rate function. Namely, using,
\begin{equation*}
\liminf_{\epsilon \rightarrow 0} \epsilon \log P\left(\sup_{0<t\leq T} \|u_{t}^{\epsilon}(y)\|_{\beta} >r\right) \geq -\inf_{u_{t}(y)\in [-r,r]^{c}} I_{1}(u_{t}(y)),
\end{equation*}
where $u_{t}(y)$ is the solution of \eqref{controlled}. Then we may deduce that for a given $\delta >0$, and sufficiently small $\epsilon >0$,
\begin{eqnarray*}\label{LDPbound}
P(\tau^{\epsilon} \leq T) &\geq& P\left(\sup_{0<t\leq T} \|u_{t}^{\epsilon}(y) \|_{\beta} >r\right)\\
&\geq& \exp \left(-\frac{1}{\epsilon} \left(\inf_{u_{t}(y) \in [-r,r]^{c}} I_{1}(u_{t}(y))+\delta\right)\right).
\end{eqnarray*}

Also, we may find,
\begin{eqnarray*}
P(\tau^{\epsilon} > 1) &=& P\left( \sup_{0< t\leq 1} \|u_{t}^{\epsilon}(y)\|_{\beta} < r\right)\\
&\leq& P\left( \sup_{0<t \leq 1} \|u_{t}^{\epsilon}(y)\|_{\beta} \leq r\right)\\
&\leq& \exp\left(-\frac{1}{\epsilon} \left( \inf_{ u_{t}^{\epsilon}(y) \in [-r,r]} I_{1}\left(u_{t}^{\epsilon}(y)\right) - \delta\right)\right).
\end{eqnarray*}
Note that by the strong Markov property of $u_{t}^{\epsilon}(y)$,
\begin{eqnarray*}
P(\tau^{\epsilon} > k+1) &=& P\left(\tau^{\epsilon} >k, \tau^{\epsilon} > k+1\right)\\
&=& \mathbb{E}_{u_{0}^{\epsilon}(y)}\left(1_{\tau^{\epsilon} >k} 1_{\tau^{\epsilon} >k+1}\right)\\
&=& \mathbb{E}_{u_{0}^{\epsilon}(y)} \left(1_{\tau^{\epsilon} > k} \mathbb{E}_{u_{0}^{\epsilon}(y)} \left(1_{\tau^{\epsilon} > k+1} |\mathcal{F}_{k}\right)\right)\\
&=& \mathbb{E}_{u_{0}^{\epsilon}(y)} \left(1_{\tau^{\epsilon} > k} P_{u_{k}^{\epsilon}(y)} (\tau^{\epsilon} > 1)\right)\\
&\leq& P(\tau^{\epsilon} >k) \sup_{u_{t}(y) \in (-r,r)} P_{u_{t}(y)} (\tau^{\epsilon} >1),
\end{eqnarray*}
so that by an inductive argument, one may deduce that for $k\in \mathbb{N}$,
\begin{equation*}
P(\tau^{\epsilon}>k) \leq \left(\sup_{u_{t}^{\epsilon}(y) \in (-r,r)} P_{u_{t}^{\epsilon}(y)}(\tau^{\epsilon}>1)\right)^{k}.
\end{equation*}

Observe that the above estimates also hold for $\tau^{\epsilon} \geq k$. Thus, as for inequality \eqref{expectationSPDE} we have,
\begin{eqnarray*}
E(\tau^{\epsilon}) \leq \sum_{k=0}^{\infty} P(\tau^{\epsilon} \geq k) &\leq& \sum_{k=0}^{\infty}\left(\sup_{u_{t}^{\epsilon}(y) \in (-r,r)} P_{u_{t}^{\epsilon}(y)}(\tau^{\epsilon}>1)\right)^{k}\\
&=& \frac{1}{1-\exp \left( -\frac{1}{\epsilon}\left( \inf_{u_{t}(y) \in [-r,r]} I_{1}\left(u_{t}(y)\right) - \delta\right)\right)},
\end{eqnarray*}
which offers another lowerbound on $P(\tau^{\epsilon} > T)$ as follows,
\begin{equation}\label{upperbound}
P(\tau^{\epsilon} > T) \leq \exp\left(-\frac{\lfloor T \rfloor}{\epsilon} \left( \inf_{ u_{t}(y) \in [-r,r]} I_{1}\left(u_{t}(y)\right) - \delta\right)\right),
\end{equation}
where $\lfloor x \rfloor$ is the greatest integer less than or equal to $x$. \\
\begin{flushright}
  $\Box$
\end{flushright}

\textbf{Proof of Theorem 2.2}: Suppose we have an attraction point in the domain, which in our case we make the point zero. Part of the definition of attraction point is the assumption that after entering its $\delta_{0}$-neighborhood for a small enough $\delta_{0}>0$, the trajectory will never leave and all trajectories will eventually converge to the attraction point as time goes to infinity. Hence for $u_{t}^{\epsilon}(y)$ starting at $y_{0}\in D \cup \partial D \setminus[-\delta_{0},\delta_{0}]$ the exit time needs to be before entering the $\delta_{0}$-neighborhood of the attraction point. Letting,
\begin{equation}\label{T1}
\tau_{1}^{\epsilon}:= \inf\{t: u_{t}^{\epsilon}(y) \in (-\delta_{0}, \delta_{0})\},
\end{equation}
 this means $\tau^{\epsilon} \leq \tau^{\epsilon}_{1}$. Using the results from the previous theorem, we observe that,
\begin{equation*}
P(\tau^{\epsilon} > T) \leq P(\tau^{\epsilon}_{1} >T) = P\left(\sup_{0<t\leq T} \|u_{t}^{\epsilon}(y)\|_{\beta} \geq \delta_{0}\right) \leq J(\delta_{0}, \epsilon, T),
\end{equation*}
from which the following lowerbound may be obtained,
\begin{equation*}
1-J(\delta_{0}, \epsilon, T) \leq P(\tau^{\epsilon} \leq T).
\end{equation*}
For $\tau^{\epsilon} \leq T$, it is required to have $T <\tau_{1}^{\epsilon}$ so that,
\begin{eqnarray*}
P(\tau^{\epsilon} \leq T <\tau_{1}^{\epsilon}) &=& P(T\geq \tau^{\epsilon}) - P(T\geq \tau_{1}^{\epsilon})\\
&\leq& J(r,\epsilon, T) + \exp\left( -\frac{1}{\epsilon}\left(\inf_{u_{t}(y)\in (-\delta_{0},\delta_{0})^{c}} I_{1}(u_{t}(y)) -\delta\right)\right),
\end{eqnarray*}
for sufficiently small $\epsilon>0$ and $\delta>0$. Using \eqref{upperbound} we may also find,
\begin{eqnarray*}
\mathbb{E}(\tau^{\epsilon}) &=& \sum_{k=0}^{\infty} P(\tau^{\epsilon} >k)\leq  \sum_{k=0}^{\infty} P(\tau_{1}^{\epsilon} >k)\\
&\leq& \sum_{k=0}^{\infty}  \exp\left(-\frac{k}{\epsilon}\left(\inf_{u_{t}(y) \in (-\delta_{0},\delta_{0})^{c}}I_{1}(u_{t}(y)) - \delta\right)\right)\\
&\leq& \frac{1}{1- \exp\left(-\frac{1}{\epsilon}\left(\inf_{u_{t}(y) \in (-\delta_{0},\delta_{0})^{c}}I_{1}(u_{t}(y)) - \delta\right)\right)}.
\end{eqnarray*}
for sufficiently small $\epsilon >0$.
\begin{flushright}
  $\Box$
\end{flushright}

We shall make the remark that in Lemma 4.2.2 and Theorem 4.4.1 of \cite{Freidlin1}, the authors assume the attraction point to be a stable equilibrium position of the domain, which means that for every neighborhood of the attraction point there is a smaller neighborhood so that if a trajectory starts in the smaller neighborhood, it will converge to the attraction point as $t \rightarrow \infty$ without leaving the larger neighborhood. In addition, they give the condition that $(b(x),n(x))<0$ for the starting point $x$ on the boundary of the domain, where $b(x)$ is the drift of the equation and $n(x)$ is the exterior normal to the boundary of the domain to ensure that the trajectories do not exit the domain. We can find the estimates in our case without requiring these assumptions. They also use the terminology of action functional denoted as $S_{0,T}(\phi)$ given to the good rate function multiplied by the speed of the large deviation principle.

\section{Exit Measure for SBM and FVP}
Here we focus on the population models, SBM and FVP, where again we have our setting in dimension one and consider $(-r,r)$ for $r>0$ as our domain of study. With regards to the bounds derived in the previous section, recall that $u_{t}^{\epsilon}(y)= \int_{0}^{y}\mu_{t}^{\epsilon}(dx)$ and $u_{t}^{\epsilon}(y)= \int_{-\infty}^{y} \mu_{t}^{\epsilon}(dx)$ for relation between the SPDE \eqref{SPDE} and SBM and FVP, respectively. Following the inner product given by \eqref{innerproduct}, we use a complete orthonormal set $\{f_{j}\}_{j}$, and sum on $j$, to deduce for each model,
\begin{eqnarray*}
P\left(\sup_{0<t\leq T} \|\mu_{t}^{\epsilon}\|_{\beta} \geq r\right) &=& P\left(\sup_{0<t\leq T} \int \sup_{y}e^{-2\beta |y|} u_{t}^{\epsilon}(y)dy \geq r^{2}\right)\\
&=& P\left(\sup_{0<t\leq T}\|u_{t}^{\epsilon}(y)\|_{\beta}  \geq r^{2} N_{1}^{-1}\right)\\
&\leq& J(r^{2}N^{-1}, \epsilon, T),
\end{eqnarray*}
where $N_{1}= \int \sup_{y} e^{-\beta |y|}dy$ and $\mu_{t}^{\epsilon}(dy)= u_{t}^{\epsilon}(y)dy$ was applied in the first step. Since $\mu^{\epsilon}_{t}$ is an empirical measure giving the size of the population up to time $t$, then the exit time for the two population models may be defined as,
\begin{equation*}
\tau^{\epsilon}:= \inf\{t: \mu_{t}^{\epsilon}((-r,r)^{c}) \geq 1\},
\end{equation*}
giving,
\begin{equation*}
P(\tau^{\epsilon} >T) = P\left(\sup_{0< t\leq T} \mu_{t}^{\epsilon}\left((-\infty, -r) \cup (r,\infty)\right) =0\right).
\end{equation*}
By the relation between SPDE \eqref{SPDE} and SBM and FVP, we may deduce that for a measurable set $B$,
\begin{equation*}
\int 1_{B}(y) du_{t}^{\epsilon}(y) = \int 1_{B}(y) \mu_{t}^{\epsilon}(dy) = \mu_{t}^{\epsilon}(B).
\end{equation*}
Similar to estimates in previous section, we may find,
\begin{eqnarray}\label{tauinequality}
P(\tau^{\epsilon} >T) &=& P\left(\sup_{0< t\leq T} \mu_{t}^{\epsilon}\left((-\infty, -r) \cup (r,\infty)\right)=0\right)\nonumber\\
&=& P\left(\sup_{0< t\leq T} \int \left(1_{(-\infty, -r)}(y) + 1_{(r,\infty)}(y)\right)du_{t}^{\epsilon}(y)=0\right)\nonumber\\
&=& 1- P\left(\sup_{0 < t\leq T} \int g(y) du_{t}^{\epsilon}(y) \geq 1\right),
\end{eqnarray}
where we have denoted,
\begin{equation*}
g(y):= 1_{(-\infty, -r)} (y) + 1_{(r,\infty)}(y).
\end{equation*}
For better presentation, let
\begin{equation*}
\tilde{P}(t,r,x):= p_{t}(r+x)-p_{t}(r-x).
\end{equation*}
Then, denoting the derivative with respect to $y$ with a prime, we have,
\begin{eqnarray*}
&&P\left( \sup_{0< t\leq T}\int g(y) du_{t}^{\epsilon}(y) \geq 1\right)\\
&=& P\left(2\sup_{0 < t\leq T}\left|\int g(y) \int_{\mathbb{R}} p_{t}'(y-x)F(x)dxdy\right|^{2}\right.\\
&&\left. +2\sup_{0 < t\leq T} \sqrt{\epsilon} \left| \int g(y) \int_{0}^{t} \int_{U} \int_{\mathbb{R}} p_{t-s}'(y-x)G(a,x,u_{s}^{\epsilon}(x))dxW(dads)dy\right|^{2} \geq 1 \right)\\
&=& P\left(2\sup_{0 < t\leq T}\left|\int_{\mathbb{R}} \tilde{P}(t,r,x) F(x)dx\right|^{2} \right.\\
&&\left.+ 2\sup_{0 < t\leq T} \sqrt{\epsilon}\left|\int_{0}^{t} \int_{U}\int_{\mathbb{R}} \tilde{P}(t-s,r,x)G(a,x,u_{s}^{\epsilon}(x))dxW(dads)\right|^{2}\geq 1\right).
\end{eqnarray*}
Analogous to bound in \eqref{Finequality} we determine,
\begin{equation*}
2\left|\int_{\mathbb{R}} \left(p_{t}(r+x) -p_{t}(r-x)\right)F(x)dx\right|^{2}\leq \frac{1}{t} \left(K_{5} + K_{6} e^{2\beta_{0}|r|}\right).
\end{equation*}
From \eqref{SBM} one may note that SBM satisfies SPDE \eqref{SPDE} with $U=\mathbb{R}, \lambda(da)=da$ and $G(a,y,u)=1_{0\leq a \leq u}+ 1_{u\leq a \leq 0}$ where $\lambda$ is the Lebesgue measure. Further notice that \eqref{FVP} implies FVP satisfies SPDE \eqref{SPDE} with $u_{t}^{\epsilon}(y)= \int_{-\infty}^{y} \mu_{t}^{\epsilon}(dx)$, $U=[0,1], \lambda(da)= da,$ and $G(a,y,u)= 1_{a<u}-u$. In both case the relation $du_{t}^{\epsilon}(y)=\mu_{t}^{\epsilon}(dy)$ along with conditions \eqref{condition1} and \eqref{condition2} hold. We proceed with,
\begin{eqnarray*}
&&P\left(\sup_{0< t\leq T} \int g(y) du_{t}^{\epsilon}(y) \geq 1\right)\\
&=& P\left(2\sqrt{\epsilon} \sup_{0<t\leq T} \left|\int_{0}^{t} \int_{\mathbb{R}} \int_{\mathbb{R}} \tilde{P}(t-s, r,x) G(a,x,u_{s}^{\epsilon}(s))dxW(dads)\right|^{2}\right. \\
 &&\left.\geq 1- \sup_{0 < t\leq T} \frac{1}{t} \left(K_{5} + K_{6} e^{2\beta_{0}|r|}\right)\right).
\end{eqnarray*}

 Applying condition \eqref{condition2}, we follow estimates as in previous section to obtain,
 \begin{eqnarray*}
 &&\sup_{0<t\leq T} \sqrt{\epsilon} \left|\int_{0}^{t} \int_{\mathbb{R}}\int_{\mathbb{R}} \tilde{P}(t-s,r,x)G(a,x,u_{s}^{\epsilon}(x))dxW(dads)\right|^{2}\\
 &\leq& \sqrt{\epsilon} \sup_{0<t\leq T}\int_{0}^{t} \int_{\mathbb{R}} \int_{\mathbb{R}} \tilde{P}(t-s,r,x)^{2} e^{3\beta_{1}|x|} dx \int_{\mathbb{R}} G(a,x,u_{s}^{\epsilon}(x))^{2} e^{-3\beta_{1}|x|} dx dads \\
 &\leq& \sqrt{\epsilon} \sup_{0<t\leq T} \int_{0}^{t} \left(\frac{K_{7}}{\sqrt{t-s}} + \frac{K_{8}}{\sqrt{t-s}} e^{3\beta_{1}|r|}\right)\int (1+u_{s}^{\epsilon}(x)^{2})e^{-3\beta_{1}|x|}dxds\\
 &\leq& K_{9}\sqrt{\epsilon} \sup_{0<s\leq T} \|u_{s}^{\epsilon}(x)\|_{\beta_{1}}^{2} \sqrt{\epsilon T} (1+e^{3\beta_{1}|r|}).
 \end{eqnarray*}
 Thus,
 \begin{eqnarray*}
 &&P\left(\sup_{0<t\leq T} \int g(y)du_{t}^{\epsilon}(y) \geq 1\right)\\
 &\leq& P\left(\sup_{0<s\leq T} \|u_{s}^{\epsilon}(x)\|_{\beta_{1}}^{2} \geq \frac{1}{K_{9}\sqrt{\epsilon T} (1+e^{3\beta_{1}|r|})}\left(1- \sup_{0<t\leq T} \frac{1}{t} \left(K_{5} + K_{6} e^{2\beta_{0}|r|}\right)\right)\right).
 \end{eqnarray*}
 Now by \eqref{probability} and \eqref{tauinequality} we arrive at \eqref{result}. For \eqref{result2} we proceed as follows.\\
 Note that $\mu_{t}^{\epsilon}$ takes values in space $\mathbb{M}_{\beta}$ or $\mathbb{P}_{\beta}$ for SBM or FVP, respectively. Let $\{f_{j}\}_{j}$ be a set of positive functions on $y$ that form a complete orthonormal system  with $f_{j}\in \mathcal{C}_{c}^{\infty}(\mathbb{R})$. Using the relation, $\mu_{t}^{\epsilon}(dy)= u_{t}^{\epsilon}(y)dy$, we obtain,
\begin{eqnarray}\label{Iequation}
\mathbb{E}<\mu_{t}^{\epsilon}(dy), f_{j}(y)>_{\beta}^{2} &=& \mathbb{E}\left(\int \sup_{y} \left(e^{-\beta |y|} f_{j}'(y)-\beta f_{j}(y) \text{sign}(y)e^{-\beta |y|}\right)u_{t}^{\epsilon}(y)dy\right)^{2}\nonumber \\
&\leq& K\mathbb{E}\left(\int \sup_{y} e^{-\beta |y|}f_{j}'(y)u_{t}^{\epsilon}(y)dy\right)^{2}\nonumber \\
&&+ K\mathbb{E}\left(\int \sup_{y}e^{-\beta |y|}f_{j}(y)u_{t}^{\epsilon}(y)dy\right)^{2} \nonumber \\
&=& I_{1} + I_{2}
\end{eqnarray}
As for $I_{1}$ we have,
\begin{eqnarray*}
I_{1}&\leq& K\left(\int_{\mathbb{R}} \sup_{y} e^{-\beta |y|} f_{j}'(y) F(y)dy\right)^{2}\\
&& + K\epsilon \mathbb{E}\left(\int_{0}^{t} \int_{\mathbb{R}} \int_{\mathbb{R}} \sup_{y} e^{-\beta |y|} f_{j}'(y) \left(1_{0\leq a\leq u_{s}^{\epsilon}(y)} + 1_{u_{s}^{\epsilon}(y)\leq a\leq 0}\right)dyW(dads)\right)^{2}\\
&&+ K\mathbb{E} \left(\frac{1}{2} \int_{0}^{t} \int_{\mathbb{R}}\sup_{y} e^{-\beta |y|} f_{j}'(y) \Delta u_{s}^{\epsilon}(y)dsdy\right)^{2}\\
&=& I_{11} + I_{12} + I_{13},
\end{eqnarray*}
where,
\begin{eqnarray*}
I_{11}= K\left(\int_{\mathbb{R}} \sup_{y} e^{-\beta |y|} f_{j}'(y) F(y)dy\right)^{2}&\leq& K\left(\int_{\mathbb{R}} \sup_{y} e^{-\beta |y|} f_{j}'(y) e^{\beta_{0}|y|} dy\right)^{2}\\
&\leq& K\left(\int_{\mathbb{R}} \sup_{y} e^{-(\beta - \beta_{0})|y|} f_{j}'(y)dy\right)^{2}\\
&\leq& K(\beta-\beta_{0})^{2} \left(\int_{\mathbb{R}} \sup_{y} e^{-(\beta-\beta_{0})|y|} f_{j}(y)dy\right)^{2}.
\end{eqnarray*}

Furthermore,
\begin{eqnarray*}
I_{12}&\leq& K\epsilon \mathbb{E} \left(\int_{0}^{t} \int_{\mathbb{R}} \left(\int_{\mathbb{R}} \sup_{y} e^{-\beta |y|} f_{j}'(y)
\left(1_{0\leq a\leq u_{s}^{\epsilon}(y)} + 1_{u_{s}^{\epsilon}(y)\leq a\leq 0}\right)dy \right)^{2}dads\right)\\
&\leq& K\epsilon \mathbb{E} \left(\int_{0}^{t} \int_{\mathbb{R}} \int_{\mathbb{R}} \sup_{y} e^{-2\beta |y|} f_{j}'(y)dy
\int_{\mathbb{R}} (1_{0\leq a\leq u_{s}^{\epsilon}(y)} + 1_{u_{s}^{\epsilon}(y)\leq a\leq 0})^{2}f'_{j}(y) dy dads\right),
\end{eqnarray*}
where using the fact that $u_{t}^{\epsilon}(y)$ is an increasing function and $\mu_{t}^{\epsilon}$ being the size of the population, we obtain,
\begin{eqnarray}\label{aestimate}
&&\int_{\mathbb{R}}\int_{\mathbb{R}} \left(1_{0\leq a \leq u_{s}^{\epsilon}(y)} + 1_{u_{s}^{\epsilon}(y) \leq a\leq 0} \right)^{2}f_{j}'(y)dyda\nonumber \\
&\leq& 2 \int_{0}^{\infty} \int_{(u_{s}^{\epsilon})^{-1}(a)}^{\infty} f_{j}'(y)dy da + 2 \int_{-\infty}^{0} \int_{-\infty}^{(u_{s}^{\epsilon})^{-1}(a)} f_{j}'(y)dyda\nonumber \\
&=& -2\int_{0}^{\infty}f_{j}((u_{s}^{\epsilon})^{-1}(a)) da + 2 \int_{-\infty}^{0} f_{j}( (u_{s}^{\epsilon})^{-1}(a)) da \nonumber \\
&=& -2 \int_{0}^{\infty} f_{j}(y) du_{s}^{\epsilon}(y) + 2 \int_{-\infty}^{0} f_{j}(y) du_{s}^{\epsilon}(y)\nonumber \\
&=& -2\int_{0}^{\infty} f_{j}(y)\mu_{s}^{\epsilon}(dy) + 2 \int_{-\infty}^{0} f_{j}(y)\mu_{s}^{\epsilon}(dy) \nonumber \\
&\leq& 2 \int_{\mathbb{R}} f_{j}(y)\mu_{s}^{\epsilon}(dy),
\end{eqnarray}
giving,
\begin{eqnarray*}
I_{12} &\leq& K\epsilon \beta \mathbb{E} \left(\int_{0}^{t} \int_{\mathbb{R}} \sup_{y}e^{-2\beta |y|} f_{j}(y)dy \int_{\mathbb{R}} f_{j}(y) \mu_{s}^{\epsilon}(dy)ds\right),
\end{eqnarray*}
where the positivity of the integrand was observed. Turning to $I_{13}$ we have,
\begin{eqnarray*}
I_{13} &=& K\mathbb{E} \left(\frac{1}{2} \int_{0}^{t} \int_{\mathbb{R}} \sup_{y} e^{-\beta |y|} f_{j}'(y) \Delta u_{s}^{\epsilon}(y)dsdy\right)^{2}\\
 &\leq& K t^{2} \mathbb{E} \left(  \int \sup_{y} \left(e^{-\beta |y|}f_{j}'(y)\right)'' \sup_{0\leq s\leq t} u_{s}^{\epsilon}(y)dy\right)^{2}\\
&=& K t^{2} \mathbb{E} \left(\int \sup_{y} e^{-\beta |y|} f_{j}'''(y) \sup_{0\leq s\leq t} u_{s}^{\epsilon}(y)dy\right)^{2}\\
&&+ 2 K\beta^{2} t^{2} \mathbb{E} \left( \int_{\mathbb{R}} \sup_{y} e^{-\beta |y|} f_{j}''(y) \sup_{0\leq s\leq t}u_{s}^{\epsilon}(y)dy\right)^{2}\\
&& + K\beta^{4} t^{2}\mathbb{E} \left(\int \sup_{y}e^{-\beta |y|} f_{j}'(y)\sup_{0\leq s\leq t}u_{s}^{\epsilon}(y)dy\right)^{2}\\
&=& I_{131} + I_{132} + I_{133}
\end{eqnarray*}

As for $I_{131}$,
\begin{eqnarray*}
I_{131} &\leq& Kt^{2}\mathbb{E} \left(\int \sup_{y} e^{-2\beta |y|}f_{j}'''(y)dy \int \sup_{0\leq s\leq t} u_{s}^{\epsilon}(y)^{2} f_{j}'''(y)dy\right)\\
&=& -K\beta^{3} \text{sign}(y) t^{2}\mathbb{E} \left(\int_{\mathbb{R}} \sup_{y} e^{-2\beta |y|}f_{j}(y)dy\int \sup_{0\leq s\leq t} \left(u_{s}^{\epsilon}(y)^{2}\right)''' f_{j}(y)dy\right).
\end{eqnarray*}

Also,
\begin{eqnarray*}
I_{132} &\leq& K\beta^{2} t^{2} \mathbb{E}\left(\int \sup_{y}e^{-\beta |y|} f_{j}''(y)\sup_{0\leq s\leq t}u_{s}^{\epsilon}(y)dy\right)^{2}\\
&\leq&  K\beta^{2} t^{2} \mathbb{E} \left(\int \sup_{y} e^{-2\beta |y|} f_{j}''(y)dy \int \sup_{0\leq s\leq t} u_{s}^{\epsilon}(y)^{2} f_{j}''(y)dy\right)\\
&\leq& K\beta^{2} t^{2} \mathbb{E}\left(\int \sup_{y} e^{-2\beta |y|} f_{j}(y)dy \int \sup_{0\leq s\leq t} \left(u_{s}^{\epsilon}(y)^{2}\right)'' f_{j}(y)dy\right),
\end{eqnarray*}
and
\begin{eqnarray*}
I_{133} &\leq & K t^{2} \mathbb{E} \left(\int \sup_{y} e^{-2\beta |y|} f_{j}'(y)dy \int \sup_{0\leq s\leq t} u_{s}^{\epsilon}(y)^{2} f_{j}'(y)dy\right)\\
&=& -K\beta t^{2} \text{sign} (y) \mathbb{E} \left( \int \sup_{y} e^{-2\beta |y|} f_{j}(y)dy \int_{\mathbb{R}} \sup_{0\leq s \leq t} \left(u_{s}^{\epsilon}(y)^{2}\right)' f_{j}(y)dy\right).
\end{eqnarray*}

Next for $I_{2}$ in \eqref{Iequation} we use the mild form and obtain,
\begin{flalign*}
I_{2}&= K\mathbb{E}\left(\sup_{y} e^{-\beta |y|}f_{j}(y)u_{t}^{\epsilon}(y)dy\right)^{2}&\\
&\leq K\left(\int \sup_{y}e^{-\beta |y|}f_{j}(y) \int p_{t}(y-x)F(x)dx dy\right)^2&\\
&+ K\mathbb{E} \left(\sqrt{\epsilon} \int \sup_{y} e^{-\beta |y|} f_{j}(y) \int_{0}^{t} \int_{U} \int_{\mathbb{R}} p_{t-s}(y-x)\left(1_{0\leq a\leq u_{s}^{\epsilon}(x)} + 1_{u_{s}^{\epsilon}(x)\leq a\leq 0}\right)dx W(dads)dy\right)^{2}&\\
&= I_{21} + I_{22},&
\end{flalign*}

where,
\begin{eqnarray*}
I_{21} &\leq& K\left(\int \sup_{y} e^{-\beta |y|} f_{j}(y) \int p_{t}(y-x)K_{1}e^{\beta_{0}|x|}dxdy\right)^{2}\\
&\leq& K \left(\int \sup_{y} e^{-\beta |y|}f_{j}(y) \frac{1}{\sqrt{t}} e^{\beta_{0}|y|}dy\right)^{2}\\
&\leq& \frac{K}{t} \left(\int_{\mathbb{R}} \sup_{y} e^{-(\beta -\beta_{0})|y|} f_{j}(y)dy\right)^{2}.
\end{eqnarray*}
Since $G(a,y,u)= 1_{0\leq a \leq u} + 1_{u\leq a\leq 0}$ satisfies conditions \eqref{condition1} and \eqref{condition2} then by applying estimate \eqref{Mlemma} along with Burkholder-Davis-Gundy inequality we arrive at,
\begin{flalign*}
I_{22}&\leq \epsilon K \mathbb{E} \left(\int_{0}^{t} \int_{U} \left(\int \sup_{y}e^{-\beta |y|} f_{j}(y) \int_{\mathbb{R}} p_{t-s}(y-x) \left(1_{0\leq a\leq u_{s}^{\epsilon}(x)} + 1_{u_{s}^{\epsilon}(x)\leq a\leq 0}\right)dxdy\right)^{2}dads\right)&\\
&\leq \epsilon K\mathbb{E} \left(\int_{0}^{t} \int_{U} \int_{\mathbb{R}}\sup_{y}e^{-\beta |y|}f_{j}(y)dy \right.&\\
&\left. \int_{\mathbb{R}}  e^{-\beta |y|} f_{j}(y) \left(\int_{\mathbb{R}} p_{t-s}(y-x) \left(1_{0\leq a\leq u_{s}^{\epsilon}(x)} + 1_{u_{s}^{\epsilon}(x)\leq a\leq 0}\right)dx\right)^{2}dydads\right)&\\
&\leq  \epsilon K \mathbb{E} \left(\int_{0}^{t} \int_{U} \int_{\mathbb{R}} \sup_{y} f_{j}(y) e^{-\beta |y|}dy
\int_{\mathbb{R}} f_{j}(y)e^{-\beta |y|} \int_{\mathbb{R}} p_{t-s}(y-x)^{2} e^{2\beta_{1} |x|}dx \right.&\\
&\left. \int_{\mathbb{R}} \left(1_{0\leq a\leq u_{s}^{\epsilon}(x)} + 1_{u_{s}^{\epsilon}(x)\leq a\leq 0}\right)^{2}e^{-2\beta_{1}|x|} dxdydads\right)&\\
&\leq \epsilon K\mathbb{E}\left(\int_{0}^{t}  \int_{\mathbb{R}} \sup_{y} e^{-\beta |y|} f_{j}(y)dy \int_{\mathbb{R}} f_{j}(y) e^{-\beta |y|} e^{2\beta_{1}|y|} \frac{1}{\sqrt{t-s}} \int_{\mathbb{R}} (1+ u_{s}^{\epsilon}(x)^{2})e^{-2\beta_{1}|x|}dxdyds\right)&\\
&\leq \epsilon KM \left(\int_{0}^{t} \frac{1}{\sqrt{t-s}} \int_{\mathbb{R}} \sup_{y} e^{-\beta |y|} f_{j}(y)dy
 \int_{\mathbb{R}} f_{j}(y) e^{-(\beta -2\beta_{1})|y|} dyds\right)&\\
&\leq \epsilon KM \sqrt{t} \left(\int_{\mathbb{R}} \sup_{y} e^{-\beta |y|} f_{j}(y)dy \int_{\mathbb{R}} f_{j}(y)e^{-(\beta-2\beta_{1})|y|} dy\right)&.
\end{flalign*}

Thus from the above bounds we can form,
\begin{flalign*}
&\mathbb{E}\left<\mu_{t}^{\epsilon}(dy), f_{j}(y)\right>_{\beta}^{2} \leq K(\beta -\beta_{0})^{2} \left(\int_{\mathbb{R}} \sup_{y} e^{-(\beta -\beta_{0})|y|}f_{j}(y)dy\right)^{2}&\\
&+K\epsilon \beta \mathbb{E} \left(\int_{0}^{t} \int \sup_{y} e^{-2\beta |y|}f_{j}(y)dy \int f_{j}(y) \mu_{s}^{\epsilon}(dy)ds\right)&\\
&- K\beta^{3} \text{sign}(y) t^{2} \mathbb{E}\left(\int \sup_{y} e^{-2\beta |y|} f_{j}(y)dy\int \sup_{0\leq s\leq t} \left(u_{s}^{\epsilon}(y)^{2}\right)'''f_{j}(y)dy\right)&\\
&+ K\beta^{2} t^{2} \mathbb{E} \left(\int \sup_{y} e^{-2\beta |y|}f_{j}(y)dy \int \sup_{0\leq s\leq t} \left(u_{s}^{\epsilon}(y)^{2}\right)''f_{j}(y)dy\right)&\\
&- K\beta t^{2} \text{sign}(y) \mathbb{E} \left(\int \sup_{y} e^{-2\beta |y|}f_{j}(y)dy \int \sup_{0\leq s\leq t} \left(u_{s}^{\epsilon}(y)^{2}\right)'f_{j}(y)dy\right)&\\
&+\frac{K}{t}\left(\int \sup_{y} e^{-(\beta -\beta_{0})|y|}f_{j}(y)dy\right)^{2}&\\
& +\epsilon K M\sqrt{t} \left(\int \sup_{y}e^{-\beta |y|} f_{j}(y)dy \int f_{j}(y) e^{-(\beta-2\beta_{1})|y|}dy\right).&
\end{flalign*}
Now summing on $j$ gives,
\begin{flalign}\label{bound}
&\mathbb{E}\|\mu_{t}^{\epsilon}\|_{\beta}^{2} \leq K\left((\beta- \beta_{0})^{2} + \frac{1}{t}\right) \int \sup_{y}e^{-2(\beta - \beta_{0})|y|}dy & \nonumber \\
&+ \epsilon KM\sqrt{t} \int \sup_{y} e^{-2(\beta -\beta_{1})|y|}dy + K\epsilon \beta \mathbb{E}\int_{0}^{t} \|\mu_{s}^{\epsilon}\|_{\beta}^{2}ds&\nonumber \\
&- K\beta^{3} \text{sign}(y) t^{2} \mathbb{E} \left(\int \sup_{y} e^{-2\beta |y|} \sup_{0\leq s\leq t} \left(u_{s}^{\epsilon}(y)^{2}\right)'''dy\right)&\nonumber \\
&+ K\beta^{2} t^{2}\mathbb{E}\left(\int_{\mathbb{R}} \sup_{y}e^{-2\beta |y|} \sup_{0\leq s\leq t} \left(u_{s}^{\epsilon}(y)^{2}\right)''dy\right)&\nonumber \\
&- K\beta t^{2} \text{sign}(y)\mathbb{E} \left(\int \sup_{y} e^{-2\beta |y|} \sup_{0\leq s\leq t} \left(u_{s}^{\epsilon}(y)^{2}\right)'dy\right).&
\end{flalign}
Notice that,
\begin{eqnarray*}
&&-K\beta^{3} \text{sign}(y) t^{2} \mathbb{E}\left(\int \sup_{y} e^{-2\beta |y|} \sup_{0\leq s\leq t} \left(u_{s}^{\epsilon}(y)^{2}\right)'''dy\right)\\
&=&-K\beta^{6} t^{2} \text{sign}(y) \mathbb{E} \int_{\mathbb{R}} \sup_{y} e^{-2\beta |y|} \sup_{0\leq s\leq t} u_{s}^{\epsilon}(y)^{2}dy\\
&=& - K\beta^{6} t^{2} \mathbb{E} \int_{0}^{\infty}\sup_{y} e^{-2\beta |y|}\sup_{0\leq s\leq t} u_{s}^{\epsilon}(y)^{2}dy\\
&&+ K\beta^{6} t^{2} \mathbb{E} \int_{-\infty}^{0} \sup_{y} e^{-2\beta |y|} \sup_{0\leq s\leq t} u_{s}^{\epsilon}(y)^{2}dy\\
&\leq& K\beta^{6} t^{2}M.
\end{eqnarray*}
After similar steps are performed for the other terms we arrive at,
\begin{eqnarray*}
\mathbb{E}\|\mu_{t}^{\epsilon}\|_{\beta}^{2} &\leq& K\left((\beta-\beta_{0})^{2} + \frac{1}{t}\right) + \epsilon KM\sqrt{t} + Kt^{2}\beta^{6}M \\
&&+ K\beta^{4}t^{2}M -K\beta^{2} t^{2}M+ K\epsilon \beta \mathbb{E} \int_{0}^{t} \|\mu_{s}^{\epsilon}\|_{\beta}^{2}ds,
\end{eqnarray*}
which by Gronwall's inequality yields,
\begin{eqnarray*}
\mathbb{E}\|\mu_{t}^{\epsilon}\|_{\beta}^{2} &\leq& K\left((\beta -\beta_{0})^{2} + \frac{1}{t}\right) + KM\epsilon \sqrt{t}+ t^{2}M\left(K\beta^{6} + K\beta^{4} -K\beta^{2}\right) \\
&& + \int_{0}^{t} \left(K\left((\beta-\beta_{0})^{2}+ \frac{1}{s}\right) + KM\epsilon \sqrt{s} + s^{2}M\left(K\beta^{6} + K\beta^{4} -K\beta^{2}\right)\right)e^{-K\epsilon \beta (t-s)}ds\\
&\leq& K\left((\beta-\beta_{0})^{2} + \frac{1}{t}\right) + KM\epsilon \sqrt{t}+ t^{2} M\left(K\beta^{6} + K\beta^{4} - K\beta^{2}\right)\\
&&+ K(\beta -\beta_{0})^{2} t + \ln t + KM\epsilon t^{\frac{3}{2}}  + t^{3} M \left(K\beta^{6} + K\beta^{4}-K\beta^{2}\right).
\end{eqnarray*}
We now proceed to the case for FVP. From \eqref{FVP} one may see that FVP satisfies SPDE \eqref{SPDE} with $U=[0,1], \lambda(da)= da,$ where $\lambda$ is the Lebesgue measure and $G(a,y,u)= 1_{a<u}-u$ obeys conditions \eqref{condition1} and \eqref{condition2}. Following the same lines of reasoning as in the previous case, we note that the same estimates for $I_{1}$ and $I_{2}$ above can be used here except for $I_{12}$ which we bound as follows.

\begin{flalign*}
I_{12}&= K\epsilon \mathbb{E} \left(\int_{\mathbb{R}} \sup_{y \in \mathbb{R}} e^{-\beta |y|} f'_{j}(y) \int_{0}^{t} \int_{0}^{1} \left(1_{a\leq u_{s}^{\epsilon}(y)}-u_{s}^{\epsilon}(y)\right)W(dads)dy\right)^{2}&\\
&\leq K\epsilon \mathbb{E} \left(\int_{0}^{t} \int_{0}^{1} \left(\int \sup_{y}e^{-\beta |y|} f_{j}'(y) \left(1_{a\leq u_{s}^{\epsilon}(y)}-u_{s}^{\epsilon}(y)\right)dy\right)^{2}dads\right)&\\
&\leq K\epsilon \mathbb{E}\left(\int_{0}^{t} \int_{0}^{1} \int_{\mathbb{R}} \sup_{y} e^{-2\beta |y|} f_{j}'(y)dy\int_{\mathbb{R}} 1_{a\leq u_{s}^{\epsilon}(y)}f_{j}'(y)dydads\right)&\\
&+ K\epsilon \mathbb{E}\left(\int_{0}^{t} \int_{0}^{1} \int_{\mathbb{R}} \sup_{y} e^{-2\beta |y|}f_{j}'(y)dy \int_{\mathbb{R}} u_{s}^{\epsilon}(y)^{2} f_{j}'(y)dydads\right)&\\
&=I_{121} + I_{122}.&
\end{flalign*}
Similar to \eqref{aestimate} we obtain,
\begin{eqnarray*}
\int_{0}^{1} \int_{\mathbb{R}} 1_{a\leq u_{s}^{\epsilon}(y)} f_{j}'(y)dyda&=& \int_{0}^{1}\int_{(u_{s}^{\epsilon})^{-1}(a)} ^{\infty} f_{j}'(y)dyda\\
&=& -\int_{0}^{1} f_{j}\left((u_{s}^{\epsilon})^{-1}(a)\right)da\\
&=& - \int_{0}^{1} f_{j}(y)\mu_{s}^{\epsilon}(dy).
\end{eqnarray*}
Therefore,
\begin{eqnarray*}
I_{121} &\leq& -K\epsilon \mathbb{E} \left(\int_{0}^{t} \int_{\mathbb{R}} \sup_{y} e^{-2\beta |y|} f_{j}'(y) dy \int_{0}^{1} f_{j}(y)\mu_{s}^{\epsilon}(dy)ds\right)\\
&=& -K\epsilon \text{sign}(y)\beta \mathbb{E} \left(\int_{0}^{t} \int_{\mathbb{R}} \sup_{y} e^{-2\beta |y|} f_{j}(y)dy \int_{0}^{1} f_{j}(y) \mu_{s}^{\epsilon}(dy)ds\right)\\
&=&-K\epsilon \beta \mathbb{E} \left(\int_{0}^{t} \int_{0}^{\infty} \sup_{y}e^{-2\beta |y|}f_{j}(y)dy\int_{0}^{1} f_{j}(y)\mu_{s}^{\epsilon}(dy)ds\right)\\
&&+ K\epsilon \beta \mathbb{E} \left( \int_{0}^{t} \int_{-\infty}^{0} \sup_{y} e^{-2\beta |y|} f_{j}(y)dy \int_{0}^{1} f_{j}(y)\mu_{s}^{\epsilon}(dy)ds\right)\\
&\leq& K \beta \epsilon \mathbb{E}\left(\int_{0}^{t} \int_{\mathbb{R}} \sup_{y} e^{-2\beta |y|}f_{j}(y)dy \int_{\mathbb{R}} f_{j}(y)\mu_{s}^{\epsilon}(dy)ds\right),
\end{eqnarray*}
again by noticing the positivity of the integrand. Moreover,
\begin{eqnarray*}
I_{122} &=& -K\beta \epsilon \text{sign}(y) \mathbb{E} \left( \int_{0}^{t} \int_{0}^{1} \int_{\mathbb{R}} \sup_{y} e^{-2\beta |y|} f_{j}(y)dy \int_{\mathbb{R}} \left(u_{s}^{\epsilon}(y)^{2}\right)' f_{j}(y)dydads\right).
\end{eqnarray*}

Now grouping these bounds yields,
\begin{flalign*}
&\mathbb{E}\left<\mu_{t}^{\epsilon}(dy), f_{j}(y)\right>_{\beta}^{2} \leq K(\beta -\beta_{0})^{2} \left(\int_{\mathbb{R}} \sup_{y} e^{-(\beta -\beta_{0})|y|}f_{j}(y)dy\right)^{2}&\\
&+K \beta \epsilon \mathbb{E} \left( \int_{0}^{t} \int_{\mathbb{R}} \sup_{y}e^{-2\beta |y|} f_{j}(y) dy\int_{\mathbb{R}} f_{j}(y)\mu_{s}^{\epsilon}(dy)ds\right)&\\
&-K \text{sign}(y) \beta \epsilon \mathbb{E} \left( \int_{0}^{t} \int_{0}^{1} \int_{\mathbb{R}} \sup_{y} e^{-2\beta |y|} f_{j}(y)dy \int_{\mathbb{R}} \left(u_{s}^{\epsilon}(y)^{2}\right)' f_{j}(y)dydads\right)&\\
&- \beta^{3} \text{sign}(y) Kt^{2} \mathbb{E}\left(\int \sup_{y} e^{-2\beta |y|} f_{j}(y)dy\int \sup_{0\leq s\leq t} \left(u_{s}^{\epsilon}(y)^{2}\right)'''f_{j}(y)dy\right)&\\
&+ K\beta^{2} t^{2} \mathbb{E} \left(\int \sup_{y} e^{-2\beta |y|}f_{j}(y)dy \int \sup_{0\leq s\leq t} \left(u_{s}^{\epsilon}(y)^{2}\right)''f_{j}(y)dy\right)&\\
&- K\beta t^{2} \text{sign}(y) \mathbb{E} \left(\int \sup_{y} e^{-2\beta |y|}f_{j}(y)dy \int \sup_{0\leq s\leq t} \left(u_{s}^{\epsilon}(y)^{2}\right)'f_{j}(y)dy\right)&\\
&+\frac{K}{t}\left(\int \sup_{y} e^{-(\beta -\beta_{0})|y|}f_{j}(y)dy\right)^{2}&\\
& +\epsilon K M\sqrt{t} \left(\int \sup_{y}e^{-\beta |y|} f_{j}(y)dy \int f_{j}(y) e^{-(\beta-2\beta_{1})|y|}dy\right),&
\end{flalign*}
which by summing on $j$ and forming bounds on terms as was performed for \eqref{bound} becomes,
\begin{eqnarray*}
\mathbb{E}\|\mu_{t}^{\epsilon}\|_{\beta}^{2} &\leq& K\left((\beta- \beta_{0})^{2} + \frac{1}{t}\right)+ \epsilon KM\sqrt{t}  +Kt^{2}\beta^{6} M + K\beta^{4}t^{2}M \\
&&- K\beta^{2} t^{2}M  +K\epsilon \beta \mathbb{E}\int_{0}^{t} \|\mu_{t}^{\epsilon}\|_{\beta}^{2}ds.
\end{eqnarray*}
Hence the Gronwall's inequality gives the same result as in SBM.
\begin{flushright}
  $\Box$
\end{flushright}

\section*{Acknowledgements}
I like to thank the anonymous referee for his/her comments that helped improve the content of the paper and also my gratitude goes to the Editor in Chief for handling the publication process.

\end{document}